\documentstyle{article}
\begin{document}
\bigskip

\begin{center}
{ \Large \bf
 On Validated Numeric of Values of Some Zeta
and L-functions, and   Applications
    }
\end{center}
\smallskip
\begin{center}
{\ Nikolaj M. Glazunov}  \end{center}
\smallskip
\begin{center}
{\rm Glushkov Institute of Cybernetics NAS \\
   03187 Ukraine Kiev-187 Glushkov prospect 40 \\
   Email:} {\it glanm@d105.icyb.kiev.ua }
 \end{center}
\smallskip
\begin{flushright}
 {\it Dedicated to Hans J. Stetter on his 70th birthday}
\end{flushright}

\section{Introduction}
 This note is for an introduction to the validated numerics
study of values of zeta and L - functions. The subject matter of
this talk lies in the area between Analysis, Algebra and Validated
Numerics~\cite{AH:IC,KKU:C}.More specifically, I want to discuss
the relations between
the algebra and analysis of zeta and L - functions and the
evaluation of values of the functions.  The talk is on some
observations and partial results related to more general cases,
and I would like to ask the permission of the listeners for
replacing the notes of my talk by a very introductory exposition.
The reason for this "going down" and taking very elementary examples
as the object for discussion lies on my side, i.e. it was useful
for the study of more general cases to review the simplest cases. \\
   Why Zeta and L-functions are considered?          \\
As noted by Don Zagier,  zeta and L-functions of various sorts are
all-pervasive objects in modern number theory and its
applications~\cite{DZ:Z}.
They appear in number theory, in algebraic geometry, in algebraic
K-theory, in representation theory and in the theory of automorphic
forms. Recently zeta and L-functions appear in the theory of
dynamical systems~\cite{P:SZ,PS}.
Many mentioned applications connected with
computations of zeta and L-functions values in "critical" points
(following Deligne~\cite{PD:V} an integer $ m $ is called critical
if neither $ m $ nor $ h-m $ is a pole of $\gamma(s) $, where
$\gamma(s) $ is a "gamma factor", which equal to an exponential
function $ A^s $ times a finite product of terms $\Gamma(1/2(s+m))$
such that the product $ L^{*}(s) = \gamma(s) L(s) $ satisfies
$ L^{*}(s) = w L^{*}(k-s) $ for some integer $ h > 0 $ and
sign $ w = \pm 1 $). Then the corresponding critical value $ L(m) $
should have the form
       $$ L(m) = A(m) \Omega (m) \; ( m \; {\rm critical}) $$

where $ \Omega (m) $ is a predictable period (the integral over some
algebraic cycle of an algebraic differential form defined over a
number field) and $ A(m) $ is an algebraic number belonging to a
predictable number field.    \\
The organization of this note is as follows: In Section 1 we remind
definitions and some properties of Riemann zeta function,
Dirichlet L-functions, Dedekind zeta function, Artin-Hasse
zeta functions and Hasse-Weil zeta function of elliptic curves.
Section 2 comments some applications of zeta and $L-$functions.
In Section 3 we consider some validated numerics aspects of
zeta and $L-$functions and interval evaluation of
Riemann zeta function for $Re(s) \geq 1$ and values of some
Dirichlet L-functions at point  $s = 1.$ Also we consider some
aspects of interval evaluation of Hasse-Weil zeta functions of
elliptic  curves at point  $s = 1.$
\section{Zeta and L - functions}
   Let us remind definitions and some properties of
Riemann zeta function, Dirichlet L-functions, Dedekind
zeta function, Artin-Hasse zeta functions and Hasse-Weil
zeta function of elliptic curves.

\subsection{Riemann's zeta function}
Let
\begin{equation}
 \sum_{n} \frac{1}{n^s} = \prod_{p} (1 - p^{-s})^{-1}, n \in
{\bf N}, \; p = 2, 3, 5, 7, \ldots (primes)
\end{equation}

be the Euler product. Euler used this formula principally as a
formal identity and principally for integer values of $s.$
Let
\begin{equation}
 \zeta(s) = \sum_{n} \frac{1}{n^s}.
\end{equation}
Dirichlet, unlike Euler, used the formulae (1.1), (1.2) with $s$ as
a real variable, and, also unlike Euler, he proved rigorously that
(1.1) is true for all real $s > 1.$
Riemann defines the function~\cite{Ti,Ed}
\begin{equation}
 \zeta(s) = \frac{\Gamma(-(s+1))}{2\pi i}
 \int\limits_{+\infty}^{+\infty} \frac{(-x)^s}{exp(x) - 1}
 \frac{dx}{x} ,
\end{equation}

where $ \int\limits_{+\infty}^{+\infty} \ldots $ is a counter
integral. The limits of integration are intended to
indicate a path of integration which begins at $+\infty,$
moves to the left down the positive real axis, circles the
origin once in the positive (counterclockwise) direction,
and returns up the positive real axis to $+\infty.$
If $\zeta(s)  $ is defined by formula (1.3), then for real values of
$s$ greater than one, $\zeta(s)  $ is equal to function (1.2).
The formula (1.3) defines a function $\zeta(s)  $ which is analytic
at all points of the complex $s-$plane except for a simple pole at
$s = 1.$ The function $\zeta(s)  $  satisfies the functional
equation
\begin{equation}
 \zeta(s)  = \Gamma(-(s+1)) (2\pi)^{s-1} 2 \sin(s\pi /2)
 \zeta(1 - s).
\end{equation}
   Let $B_k$ be the $k$-th Bernoulli number in even numeration
$B_0 = 1, B_1 = -\frac{1}{2}, \ldots. $ For $n \in {\bf Z}$
\quad Euler have found that
$$ \zeta(2n) = \frac{(2\pi)^{2n}(-1)^{n+1}B_{2n}}{2(2n)!}.$$
If $n > 0$ then $\zeta(-2n) = 0.$ Functional equation (1.4)
gives for $m \in {\bf N}$ that
$$ \zeta(1-2m) = (-1)^{1-2m}\frac{B_{2m}}{2m}.$$

\subsection{Dirichlet L-functions}
  Let $p$ be a prime. Let $\omega^{p-1} = 1.$ Let $g$
be a primitive root $\pmod p.$ Let $\nu(n)$ be the index of
natural $n$ relatively to the primitive root $g$, i.e.
$g^{\nu(n)} \equiv n \pmod p.$   \\
Elementary Dirichlet character to a $\pmod p$ is the function
$$ \chi (n) = \omega^{\nu(n)}, \; \chi (n) = 0 \; \mbox{if}
\; p \mid n. $$
There are $\varphi(p) = p - 1$ distinct functions:
$\chi (n) = e^{\frac{2\pi im}{p-1}}, \; 1 \leq m \leq p - 1.$
Let $s$ be a positive real. Let
\begin{equation}
 L_{\omega}(s) =  \sum_{n=1}^{\infty} \chi (n)n^{-s}.
\end{equation}
The function $L_{\omega}(s)$ is called a Dirichlet L-function. \\
Let $\chi(n)$ be a Dirichlet character to a $\pmod q$.
In the case there exists a Dirichlet L-function
$L(s,\chi)$~\cite{Da}.
Let  $e_{q}(m) = e^\frac{2\pi im}{q}.$ Let
$$ \tau(\chi) = \sum_{m=1}^{q} \chi(m) e_{q}(m) $$
be a Gauss sum. Let now $\chi $ be a primitive character, and
let~\cite{Da} $$ \alpha = \left\{ \begin{array}{ll}
                        0, & \mbox{if $\chi(-1) = 1 $} \\
                        1, & \mbox{if $\chi(-1) = -1 $}
                        \end{array}
                \right. $$
Let
$$ \xi(s,\chi) = (\frac{\pi}{q})^{(-1/2)(s + \alpha)}
\Gamma(\frac{1}{2}(s + \alpha))L(s,\chi).$$
Then
\begin{equation}
     \xi(1 - s,\overline{\chi}) =
 \frac{i^{\alpha}q^{1/2}}{\tau(\chi)}\xi(s,\chi)
\end{equation}
is the functional equation for $L(s,\chi).$
Dirichlet $L-$functions are a partial case of Dirichlet series.

\subsection{Dedekind zeta functions}
  Dedekind extended the Riemann zeta function from natural
numbers $n$ to all integer divisors (ideals) of an algebraic
number field $K.$ The Dedekind zeta function $\zeta_{K}(s)$
is defined by the formula
\begin{equation}
   \zeta_{K}(s) = \sum_{a} \frac{1}{N(a)^s}.
\end{equation}
Here $a$ runs over all integer divisors (ideals) of the field
$K$ and $N(a)$ is the norm of $a.$ \\
Let $m \in {\bf N}, \; \xi^m = 1.$ A field $K$ is the $m-$cyclotomic
field if $K = {\bf Q}(\xi), \; \xi^m = 1.$
Let $m = p,$ where $p$ is a prime. In the case
$$ \zeta_{K}(s) = G(s) \prod_{\chi} (1 - \frac{\chi(p)}{p^s})^{-1}.$$
Here $G(s) = \prod_{\rho \mid p} (1 -
\frac{1}{N(\rho)^s})^{-1}$ and $\rho$ runs over all prime
divisors of the field ${\bf Q}(\sqrt[p]{1}), \; \rho \mid p.$
Let $\chi$ be all Dirichlet characters $\pmod p.$     \\
{\bf Proposition.}~\cite{BSh} Suppose $K$ is the
$p-$cyclotomic field. Then the Dedekind's zeta function of $K$
is the product of $G(s)$ and $\varphi(p)$ Dirichlet
$L-$functions runs over all Dirichlet characters $\chi:$
$$ \zeta_{K}(s) = G(s) \prod_{\chi} L(s,\chi). $$   \\
Let $p$ be a positive integer, $ p \equiv 1 \pmod 4$, and
$ K = {\bf Q}(\sqrt{p}).$ From Siegel's~\cite{Si} follows,
that
$$ 2\zeta_{K}(-1) = \frac{1}{15}
\sum_{1 \leq b < \sqrt{p}, \, b \, odd} \sigma_1(\frac{p - b^2}{4}),$$
where $\sigma_1(n) $ is the sum of the divisors of $n.$

\subsection{Zeta  and L-functions of elliptic curves}~\cite{Ta:AE}
 Let $E/{\bf Q}$ be an elliptic curve given in Weierstrass form
by an equation
\begin{equation}
\label{WN}
E: y^2+a_1 xy + a_3 y = x^3 + a_2 x^2 + a_4 x + a_6 ,
\end{equation}
and let $b_{2} , b_4, b_6, b_8, c_4, c_6, \Delta , j$ be the
usual associated quantities~\cite{Ta:AE}. Let now~(\ref{WN})be a
global minimal Weierstrass equation for $E$ over {\bf Z.} For each
prime $p$ the reduction~(\ref{WN}) $\pmod p$ defines a curve $E_p$ over
the prime field ${\bf F}_p$. Let $A_p$ denote the number of
 points of $E_p$ rational over ${\bf F}_p$. Let
$$ t_p = 1 + p - A_p. $$
If $p \not| \Delta,$ then $t_p$ is the trace of Frobenius
and satisfies $|T_p| \leq 2\sqrt{p}.$ In the case
{\it Artin-Hasse zeta function} of the elliptic curve $E_p$
is:
\begin{equation}
\label{ArH}
\zeta_{E_p}(s) =
\frac{1 - t_{p}p^{-s} + p^{1-2s}}{(1 - p^{-s})(1 - p^{1 - s})}.
\end{equation}
 If $p \mid \Delta,$ then $E_p$ is not an elliptic curve and has
a singularity $S.$  In the case
$$ t_p = \left\{ \begin{array}{ll}
                        0, & \mbox{if $S $ is a cusp,} \\
                        1, & \mbox{if $S $ is a node,} \\
                       -1, & \mbox{if $S $ is a node with tangent quadratic
                over ${\bf F}_p$}
                        \end{array}
                \right. $$
The {\it Hasse-Weil } $L-$function of $E/{\bf Q}$ is defined by equation
\begin{equation}
\label{HW}
\L_{E}(s) = \prod_{p \, | \Delta}
\frac{1 }{(1 - t_p p^{-s})} \prod_{p \, \not | \, \Delta}
\frac{1}{1 - t_{p}p^{-s} + p^{1-2s}}
\end{equation}
.

\section{Evaluation of values of zeta and L - functions and
their applications}
 Many results in number theory are based on evaluation of values
of Riemann zeta and Dirichlet L-functions (zeros of zeta and
L-functions, prime number theorem, distribution of primes in
arithmetic progressions, density theorems). Here we remind
some another applications (please, see also very interesting
paper~\cite{DZ:Z}).  \\

{\bf Evaluation of the index of elliptic operator and
zeta functions}~\cite{Se:P,At:GA}       \\

Let $A$ be a self adjoint positive pseudo-differential
(Calderon-Zygmund) operator. Let $A$ has no eigenvalues $\lambda$
on negative real axis. We can define
$$ A^s = \frac{1}{2\pi}\int_{\Gamma} \lambda^s (A - \lambda)^{-1}
d\lambda, $$
where $\Gamma$ is the counter which begin at $-\infty$, moves to
the right up the negative real axis, circles the origin once
in the clockwise direction, and returns down the negative real axis
to $-\infty.$ \\
 Let $A = -\frac{d^2}{dx^2} + P$ be the
operator on unit circle $S^1,$ $Pf = \frac{1}{2\pi}\int_{0}^{2\pi}
f(x)dx. $ Eigenvalues of $A$ are natural numbers $1,1,1,4,4, \ldots,
m^2, m^2, \ldots$ and if $\lambda_j$ is the $j$ eigenvalue of $A,$
then $tr A^s = Z(s) = 1 + 2\sum_{1}^{\infty} \lambda_j^s =
1 + 2\zeta(-2s).$
Suppose now that $D$ is an elliptic operator and $D^*$ is the adjoint
of $D.$ Then
$$ \Delta_0 = 1 +D^* D, \; \Delta_1 = 1 + DD^* $$
are positive self-adjoint operators. Let
$$ \zeta_i (s) = tr \Delta_{i}^{-s}, i = 0,1. $$
Then the difference of $\zeta_0$ and $\zeta_1$ gives the index of the
operator $D.$         \\

{\bf The symplectic volume of moduli space}~\cite{Do,Wi}       \\

  Let ${\cal M}(n,d)$ be the moduli space of (semi)stable
holomorphic vector bundles of rank $n$, degree $d$ and fixed
determinant on a compact Riemann surface $\Sigma.$
The symplectic volume of ${\cal M}(2,1)$ is given by
$$ vol({\cal M}(2,1)) = (1 - \frac{1}{2^{2g-3}})
\frac{\zeta(2g-2)}{2^{g-2}\pi^{2g-2}}.$$
  Here $g$ is the genus of $\Sigma.$ If $n, d, g $ are sufficiently
large, then an evaluation of $ vol({\cal M}(n,d)) $ is desirable. \\

{\bf Generalized Riemann conjecture }   \\

   Let $s = \sigma + it.$ Let $L(s,\chi)$ be a Dirichlet
$L-$function. Then
\begin{equation}
\label{GRC}
 L(s,\chi) \ne 0 \; \mbox{if} \; \sigma > \frac{1}{2}.
\end{equation}

{\bf The Hilbert modular group and the volume of
orbifolds}~\cite{H:H}          \\

   Let $K$ be a real quadratic field over $\bf Q$ and $R$
the ring of algebraic integers in $K.$ Let $H$ be the upper half
plane of $\bf C.$ The group $SL(2,R)$ acts on $H \times H.$
The Hilbert modular group $\Gamma = SL(2,R)/(1,-1) $ acts
on $H \times H $ effectively. The volume of
$ X = \Gamma \backslash (H \times H) $ is given by
$$ vol(X) = 2\zeta_{K}(-1). $$

{\bf The Birch and Swinnerton-Dyer conjecture for elliptic
curves}~\cite{Ta:AE,CW,BGZ}          \\

The Hasse-Weil $L-$function (\ref{HW}) converges for $Re(s) >
\frac{3}{2}$ because it is dominated by the product
$(\zeta(s - \frac{1}{2}))^2.$ \\
The Birch and Swinnerton-Dyer conjecture for elliptic
curves is \\
{\bf BSDC1}  Let $r$ be the rang of $E(\bf Q).$ Then
$L(E,s)$ has a zero of order $r$ at $s = 1.$ In particular,
$L(E,s)$ has the Taylor expansion
$$ L(E,s) = C_E(s-1)^r + \; \mbox{higher terms with some
constant} \; C_E \neq 0.  $$
The second Birch and Swinnerton-Dyer conjecture for elliptic
curves concerns this constant.
Arithmetic investigations involve derivatives of $L(E,s)$
and computation of their values~\cite{GZ}.

\section{Validated numerics of values of zeta and L - functions}
\subsection{Validated numerics and zeta functions}
  Remind at first the interval interpretation of inequality
$ f(z) \neq 0.$ It appears in formula (\ref{GRC}) and in many
another cases in the theory of zeta and $L-$functions.  The
inequality has the interpretation: \\
(i) If $ f(z) \in {\bf R}$ then $(f(z) > 0) \vee (f(z) < 0)$
and $f(x) \in A$ for a real interval $A, 0 \not \in A.$  \\
(ii) If $ f(z) \in {\bf C}, \; f(z) \not \in {\bf R}$ and
we use complex rectangular arithmetic $I{\bf C}$~\cite{AH:IC}, then
$(Re(z) > 0) \vee (Re(z) < 0) \vee (Im(z) > 0) \vee (Im(z) < 0).$
So $f(z) \in A, \; A \in I{\bf C}, \; 0 \not \in A.$
  Side by side with rectangular arithmetic there are complex circular
arithmetic~\cite{AH:IC} and complex sector arithmetic~\cite{KU}.
In some cases they are very useful for evaluation of
$ f(z) \neq 0.$ But as this requires the consideration of numerous
different cases we omit it.
 In preceding sections we  have reminded some results about
values of zeta functions at integer (critical) points.
 A common problem of validated numerics to find optimal interval
evaluation of the value of a function with interval variables and
parameters. In many cases for the purpose we can use:
\begin{verse}
   (i) programming languages and compilers for scientific
computations (Pascal-XSC, ACRITH-XSC, Fortran-XSC, Oberon-XSC); \\
  (ii) computer algebra systems with interval packages;         \\
 (iii) highly accurate, extended computer arithmetics with
standard and special functions for real and interval
data~\cite{KKU:C}.    \\
\end{verse}
 H. Stetter gave in~\cite{St:VC} and in his another papers
applications of computer algebra in interval methods.
In the frame of zeta and $L-$functions we can apply in some
cases algebraic-analytic approach and corresponding computer
algebra for exact presentation of zeta or $L-$function as
a finite sum and a remaining term. Then the sum and the
remaining term are interval evaluated~\cite{G:IA}.
Computer algebra systems include now in their standard library
functions some kinds of Zeta functions. For instance
Maple V Release 4 includes JacobiTheta1-Theta3 (Theta functions
are connected with Zeta functions), WeierstrassZeta, Riemann Zeta
and Hurvitz Zeta. So we can in most interesting cases compute
(not guaranteed) values of the functions. In some cases  the
Maple computing $ Zeta(n,a,m), n \in {\bf Z}, a, m \in {\bf Q} $ in
closed form (by $ \pi, a, \ln $). Standard functions and
packages of Maple give possibility to compute values of some
Zeta and L-functions that do not include in standard Maple
library. In any case we have to interval evaluate the expression
for zeta or $L-$function or the closed form of value of
zeta or $L-$function.

\subsection{Interval evaluation}
\subsubsection{Riemann zeta}

 Consider the series  $\zeta(s) = \sum_{n} \frac{1}{n^s}$
for complex values of $s$ with $Re(s) \geq 1.$ Interval
evaluation of the $\zeta(s)$ can be taken from the result
of Backlund (a little bit reformulated). The result is
based of Euler-Maclaurin summation.\\
{\bf Proposition} (Backlund) Let $N$ be natural $> 1.$
Let $s = \sigma + it$ and let $\sigma \geq 1.$ Let
$B_{2k}$ be the Bernoulli numbers in even numeration, \\
$ S(N-1,s) = \sum_{n=1}^{N-1} n^{-s} + \frac{N^{1-s}}{s-1},
\\
 B(N,k,s)
  = \frac{1}{2}N^{-s} + \frac{B_2}{2}sN^{-s-1} + \cdots
+ \frac{B_{2k}}{(2k)!}s(s+1) \ldots (s+2k-2)N^{-N-2k+1} .$ \\
Then
\begin{equation}
 \label{EMS}
\zeta(s) = S(N-1,s) + B(N,k,s) + R_{2k},
\end{equation}
where $$|R_{2k-2}| \leq \left| \frac{s+2k-1}{\sigma + 2k -1}
\right| |B_{2k} \; \mbox{term of }(\ref{EMS})|. $$

\subsubsection{Dirichlet series}

Let $K = {\bf Q}(\sqrt{D})$ be a real quadratic field with
positive integer squarefree $D$ and $\chi(n) =
\left( \frac{\Delta}{n} \right)$ be the Kronecker symbol. Here
$$\Delta = \left\{  \begin{array}{ll}
                        D, & \mbox{ $ D \equiv 1 \pmod 4 $ ,} \\
                        4D, & \mbox{ $ D \equiv 2,3 \pmod 4  $.}
                        \end{array}
                \right. $$
Let $A = \pi/\Delta, \; E(x) = \int_{x}^{\infty} \frac{e^{-t}}{t}dt,
\; erfc(x) = \frac{2}{\pi} \int_{x}^{\infty} e^{-t^2}dt. $
  \\
{\bf Proposition}~\cite{WB,PB}. \\
$$ L(1,\chi) = \frac{1}{\sqrt{\Delta}} \sum_{n=1}^{m} \chi(n) E(An^2)
\; + \sum_{n=1}^{m} \left( \frac{\chi(n)}{n} \right)
erfc(n \sqrt{A}) + R_m,$$
where $|R_m| < \frac{\Delta^{3/2}}{\pi^2} \frac{e^{-Am^2}}{m^3}.$

\subsubsection{Remarks about $L-$functions of elliptic curves.}

   From the work of A. Wiles, R. Taylor and A. Wiles, and work
of F. Diamond it is known that (semistable) elliptic curves
over $E/{\bf Q}$  are modular. Knowing the modularity of
$E/{\bf Q}$  is equivalent to the existence of a modular form $f$
on $\Gamma_{0}(N)$ for some natural value $N,$ which we write
$f = \sum a_n q^n.$ The $L-$function of $E$ is thus given by
the Mellin transform of $f, \; L(f,s) = \sum a_n q^n.$ In 
particular, the behavior of $L(E,s)$ at $s = 1$ can be deduced
from modular properties of $E.$ \\
Let $E/{\bf Q}$ be a modular elliptic curve and the global minimal
model of the $E/{\bf Q}$ has prime conductor $l$. Let $p$ be a prime
and $A_p$ (as in (1.4)) be the number of points of $E_p$ in 
${\bf F}_p.$ Then there is exists a modular form $f$ on
$\Gamma_{0}(l), \;  f = \sum_{n=1}^{\infty} a_n q^n,$ where
$a_p, \; p \neq l,$ equals $p + 1 - A_p.$ Under these assumptions 
it seems that results of~\cite{BGZ} gives expressions for 
$L(E,1)$ and $L^{'}(E,s),$ which can evaluate by validated numerics
methods.

\end{document}